\documentclass[12pt,a4paper]{article}
\usepackage{amsmath,amsthm,amsfonts,amssymb,bbm}
\usepackage{graphicx,psfrag,subfigure,epsfig}
\usepackage{cite}
\usepackage{hyperref}

\title{Random Growth Models}

\author{Patrik L. Ferrari\thanks{Institute for Applied Mathematics, Bonn University, Endenicher Allee 60,\newline 53115 Bonn, Germany; E-mail:~\texttt{ferrari@uni-bonn.de}} and Herbert Spohn\thanks{Department of Mathematics, Technical University of Munich,\newline
Boltzmannstr. 3, 85748 Garching, Germany; E-mail:~\texttt{spohn@ma.tum.de}}}

\date{13.3.2010}

\begin{document}
\maketitle
\sloppy

\begin{abstract}
The link between a particular class of growth processes and random
matrices was established in the now famous 1999 article of Baik,
Deift, and Johansson on the length of the longest increasing
subsequence of a random permutation~\cite{BDJ99}. During the past ten years, this
connection has been worked out in detail and led to an improved
understanding of the large scale properties of one-dimensional
growth models. The reader will find a commented list of references
at the end. Our objective is to provide an introduction highlighting
random matrices. From the outset it should be emphasized that this
connection is fragile. Only certain aspects, and only for specific
models, the growth process can be reexpressed in terms of partition
functions also appearing in random matrix theory.
\end{abstract}

\newpage
\section{Growth models}\label{sect1} 
\begin{figure}
\begin{center}
\epsfig{file=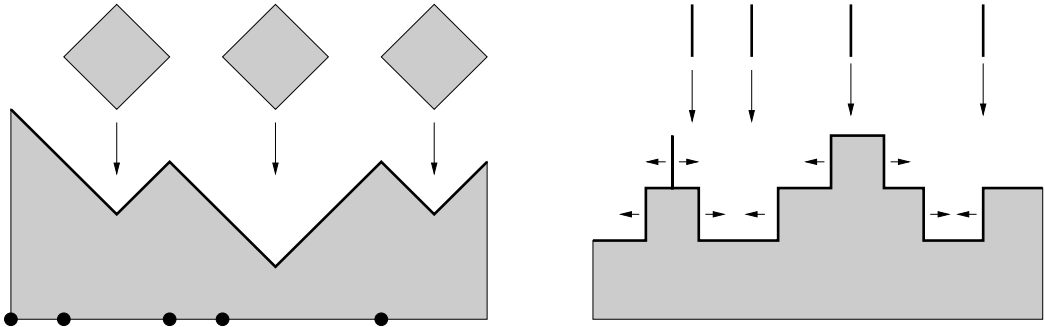,width=\textwidth}
\caption{Growth of TASEP interface (left) and continuous time PNG (right)}
\label{FigGrowth}
\end{center}
\end{figure}

A growth model is a stochastic
evolution for a height function $h(x,t)$, $x$ space, $t$ time. For
the one-dimensional models considered here, either $x\in \mathbb{R}$
or $ x\in\mathbb{Z}$. We first define the TASEP\footnote{Here we use
the height function representation. The standard particle
representation consists in placing a particle at $x$ if $h(x+1)-h(x)=-1$
and leaving empty if \mbox{$h(x+1)-h(x)=1$}.} (totally asymmetric simple
exclusion process) with parallel updating, for which
$h,x,t\in\mathbb{Z}$. An admissible height function, $h$, has to
satisfy $h(x+1)-h(x)=\pm 1$. Given $h(x,t)$ and $x_\ast$ a local
minimum of $h(x,t)$, one defines
\begin{equation}\label{2.1}
h(x_\ast, t+1)=
\left\{
\begin{array}{ll}
h(x_\ast,t)+2 & \hbox{with probability } 1-q,\; 0\leq q\leq 1\,, \\
h(x_\ast,t) & \hbox{with probability } q
\end{array}
\right.
\end{equation}
independently for all local minima, and $h(x,t+1)=h(x,t)$ otherwise, see Figure~\ref{FigGrowth} (left).
Note that if $h(\cdot,t)$ is admissible, so is $h(\cdot,t+1)$.

There are two limiting cases of interest. At $x_\ast$ the waiting
time for an increase by $2$ has the geometric distribution $(1-q)
q^n$, $n=0,1,\ldots$. Taking the $q\to 1$ limit and setting the time
unit to $1-q$ one obtains the exponential distribution of mean $1$.
This is the time-continuous TASEP for which, counting from the
moment of the first appearance, the heights at local minima are increased
independently by $2$ after an exponentially distributed waiting
time. Thus $h,x\in\mathbb{Z}$ and $t\in\mathbb{R}$.

The second case is the limit of rare events (Poisson points), where
the unit of space and time is $\sqrt{q}$ and one takes $q\to 0$.
Then the lattice spacing and time become continuous. This limit,
after a slightly different height representation (see
Section~\ref{sect5} for more insights), results in the polynuclear
growth model (PNG) for which $h\in\mathbb{Z}$, $x,t\in\mathbb{R}$.
An admissible height function is piecewise constant with jump size
$\pm 1$, where an increase by $1$ is called an up-step and a decrease by
$1$ a down-step. The dynamics is constructed from a
space-time Poisson process of intensity $2$ of nucleation events.
$h(x,t)$ evolves deterministically through (1) up-steps move to the
left with velocity $-1$, down-steps move to the right with velocity
$+1$, (2) steps disappear upon coalescence, and (3) at points of the space-time
Poisson process the height is increased by $1$, thereby nucleating
an adjacent pair of up-step and down-step. They then move
symmetrically apart by the mechanism described under (1), see Figure~\ref{FigGrowth} (right).

The TASEP and PNG have to be supplemented by initial conditions and,
possibly by boundary conditions. For the former one roughly divides
between macroscopically flat and curved. For the TASEP examples
would be $h(x)=0$ for $x$ even and $h(x)=1$ for $x$ odd, which has
slope zero, and $h(x+1)-h(x)$ independent Bernoulli random variables
with mean $m$, $|m|\leq 1$, which has slope $m$. An example for a
curved initial condition is $h(0)=0$, $h(x+1)-h(x)=-1$ for
$x=-1,-2,\ldots$, and $h(x+1)-h(x)=1$ for $x=0,1,\ldots$.

What are the quantities of interest? The most basic one is the
macroscopic shape, which corresponds to a law of large numbers for
\begin{equation}\label{2.2}
\frac{1}{t}h\big([yt],[st]\big)
\end{equation}
in the limit $t\to\infty$ with $y\in \mathbb{R}$, $s\in\mathbb{R}$,
and $[\cdot]$ denoting integer part. From a statistical mechanics
point of view the shape fluctuations are of prime concern. For
example, in the flat case the surface stays macroscopically flat and
advances with constant velocity $v$. One then would like to
understand the large scale limit of the fluctuations
$\{h(x,t)-vt,(x,t)\in\mathbb{Z}\times\mathbb{Z}\}$. As will be
discussed below, the excitement is triggered through non-classical scaling exponents and non-Gaussian limits.

In 1986 in a seminal paper Kardar, Parisi, and Zhang (KPZ) proposed
the stochastic evolution equation~\cite{KPZ86}
\begin{equation}\label{2.3}
\frac{\partial}{\partial t}h(x,t)=\lambda\Big(\frac{\partial}
{\partial x}h(x,t)\Big)^2 +\nu_0 \frac{\partial^2}{\partial
x^2}h(x,t)+\eta (x,t)
\end{equation}
for which $h,x,t\in\mathbb{R}$. $\eta(x,t)$ is space-time white
noise which models the deposition mechanism in a moving frame of
reference. The nonlinearity reflects the slope dependent growth
velocity and the Laplacian with $\nu_0>0$ is a smoothing
mechanism. To make the equation well-defined one has to introduce
either a suitable spatial discretization or a noise covariance
\mbox{$\langle \eta(x,t)\eta(x',t')\rangle=g(x-x') \delta(t-t')$} with
$g(x)=g(-x)$ and supported close to~0. KPZ argue that, according to
(\ref{2.3}) with initial conditions $h(x,0)=0$, the surface width
increases as $t^{1/3}$, while lateral correlations increase as
$t^{2/3}$. It is only through the connection to random matrix theory
that universal probability distributions and scaling functions have
become accessible.

Following TASEP, PNG, and KPZ as guiding examples it is easy to
construct variations. For example, for the TASEP one could introduce
evaporation through which heights at local maxima are decreased by 2. The
deposition could be made to depend on neighboring heights. Also
generalizations to higher dimensions, $x\in\mathbb{R}^d$ or
$x\in\mathbb{Z}^d$, are easily accomplished. For the KPZ equation
the nonlinearity then reads $\lambda\big(\nabla_x h(x,t)\big)^2$ and
the smoothening is $\nu_0 \Delta h(x,t)$. For the PNG model in
$d=2$, at a nucleation event one generates on the existing layer a
new layer of height one consisting of a disk expanding at unit
speed. None of these models seem to be directly connected to random matrices.

\section{How do random matrices appear?}\label{sect2} Let us consider
the PNG model with the initial condition $h(x,0)=0$ under the
constraint that there are no nucleations outside the interval
$[-t,t]$, $t\geq 0$, which also refered to as PNG droplet, since the
typical shape for large times is a semicircle. We study the probability distribution
$\mathbb{P}(\{h(0,t)\leq n\})$, which depends only on the nucleation
events in the quadrant \mbox{$\{(x,s)|\;|x|\leq s, s-t\leq x\leq
t-s\,,\ 0\leq s\leq t\}$}. Let us denote by
\mbox{$\omega=(\omega^{(1)},\ldots,\omega^{(n)})$} a set of nucleation
events and by $h(0,t;\omega)$ the corresponding height. The order of
the coordinates of the $\omega^{(j)}$'s in the frame $\{x=\pm t\}$
naturally defines a permutation of $n$ elements. It can be seen that
$h(0,t;\omega)$ is simply the length of the longest increasing
subsequence of that permutation. By the Poisson statistics of
$\omega$, the permutations are random and their length is Poisson
distributed. By this reasoning, somewhat unexpectedly, one finds
that $\mathbb{P}\big(\{h(0,t)\leq n\}\big)$ can be written as a matrix
integral~\cite{BR01a}. Let $\mathcal{U}_n$ be the set of all unitaries on
$\mathbb{C}^n$ and $\mathrm{d}U$ be the corresponding Haar measure.
Then
\begin{equation}\label{3.1}
\mathbb{P}\big(\{h(0,t)\leq n\}\big)= e^{-t^2} \int_{\mathcal{U}_n} \mathrm{d}U
\exp [t\, \mathrm{tr} (U+U^\ast)]\,.
\end{equation}

(\ref{3.1}) can also be expressed as Fredholm determinant on
$\ell_2=\ell_2(\mathbb{Z})$. On $\ell_2$ we define the linear
operator $B$ through
\begin{equation}\label{3.2}
(Bf)(x)=-f(x+1)-f(x-1)+\frac{x}{t} f(x)
\end{equation}
and denote by $P_{\leq 0}$ the spectral projection onto $B\leq 0$.
Setting $\theta_n(x)=1$ for $x>n$ and $\theta_n(x)=0$ for $x\leq n$,
one has
\begin{equation}\label{3.3}
\mathbb{P}\big(\{h(0,t)\leq n\}\big)= \det (\mathbf{1}-\theta_n P_{\leq 0}
\theta_n)\,.
\end{equation}
Such an expression is familiar from GUE random matrices. Let
$\lambda_1<\ldots<\lambda_N$ be the eigenvalues of an $N\times N$
GUE distributed random matrix. Then
\begin{equation}\label{3.4}
\mathbb{P}\big(\{\lambda_N\leq y\}\big)= \det (\mathbf{1}-\theta_y P_N \theta_y)\,.
\end{equation}
Now the determinant is over $L^2(\mathbb{R},\mathrm{d}y)$. If one
sets
$H=-\frac{1}{2}\frac{\mathrm{d}^2}{\mathrm{d}y^2}+\frac{1}{2N}y^2$,
then $P_N$ projects onto $H\leq \sqrt{N}$.

For large $N$ one has the asymptotics
\begin{equation}\label{3.5}
\lambda_N\cong 2N +N^{1/3} \xi_2
\end{equation}
with $\xi_2$ a Tracy-Widom distributed random variable, that is,
\begin{equation}\label{3.5a}
\mathbb{P}(\{\xi_2\leq s\})=F_2(s):=\det(\mathbf{1}-\chi_s K_{\rm Ai} \chi_s),
\end{equation}
with $\det$ the Fredholm determinant on $L^2(\mathbb{R},\mathrm{d}x)$, $\chi_s(x)=\mathbf{1}(x>s)$, and $K_{\rm Ai}$ is the Airy kernel (see~(\ref{4.4}) below).
Hence it is not so surprising that for the height of the PNG model one obtains~\cite{PS00b}
\begin{equation}\label{3.6}
h(0,t)= 2t + t^{1/3} \xi_2
\end{equation}
in the limit $t\to\infty$. In particular, the surface width
increases as $t^{1/3}$ in accordance with the KPZ prediction. The
law in (\ref{3.6}) is expected to be universal and to hold whenever
the macroscopic profile at the reference point, $x=0$ in our
example, is curved. Indeed for the PNG droplet \mbox{$h(x,t)\cong
2t\sqrt{1-(x/t)^2}$} for large $t$, $|x|\leq t$.

One may wonder whether (\ref{3.6}) should be regarded as an accident or whether there is a deeper
reason. In the latter case, further height statistics might also be representable through matrix integrals.

\section{Multi-matrix models and line ensembles}\label{sect3}
For curved initial data the link to random matrix theory can be
understood from underlying line ensembles. They differ from case to
case but have the property that the top line has the same statistics
in the scaling limit. We first turn to random matrices by
introducing matrix-valued diffusion processes.

Let $B(t)$ be GUE Brownian motion, to say $B(t)$ is an $N\times N$
hermitian matrix such that, for every $f\in \mathbb{C}^N$, $t\mapsto
\langle f,B(t) f\rangle$ is standard Brownian motion with variance
$\langle f,f\rangle^2\min (t,t')$ and such that for every unitary $U$
it holds
\begin{equation}\label{4.1}
U B(t) U^\ast =B(t)
\end{equation}
in distribution. The $N\times N$ matrix-valued diffusion, $A(t)$, is
defined through the stochastic differential equation
\begin{equation}\label{4.2}
\mathrm{d} A(t)=-V'(A(t)) \mathrm{d}t + \mathrm{d}B (t)\,,\quad A(0)=A,
\end{equation}
with potential $V:\mathbb{R}\to\mathbb{R}$. We assume $A=A^\ast$,
then also $A(t)=A(t)^\ast$ for all $t\geq 0$. It can be shown that
eigenvalues of $A(t)$ do not cross with probability 1 and we order
them as $\lambda_1(t)<\ldots<\lambda_N(t)$. $t\mapsto
\big(\lambda_1(t),\ldots,\lambda_N(t)\big)$ is the line ensemble
associated to (\ref{4.2}).

In our context the largest eigenvalue, $\lambda_N(t)$ is of most
interest. For \mbox{$N\to\infty$} its statistics is expected to be independent of
the choice of $V$. We first define the limit process, the Airy$_2$
process $\mathcal{A}_2(t)$, through its finite dimensional
distributions. Let
\begin{equation}\label{4.3}
H_{\mathrm{Ai}}=-\frac{\mathrm{d}^2}{\mathrm{d}y^2}+y
\end{equation}
as a self-adjoint operator on $L^2(\mathbb{R},\mathrm{d}y)$. The
Airy operator has $\mathbb{R}$ as spectrum with the Airy function Ai
as generalized eigenfunctions, \mbox{$H_{\mathrm{Ai}}
\mathrm{Ai}(y-\lambda)=\lambda\mathrm{Ai}(y-\lambda)$}.
In particular the projection onto $\{H_{\mathrm{Ai}}\leq 0\}$ is given by the Airy kernel
\begin{equation}\label{4.4}
K_{\mathrm{Ai}}(y,y')=\int^\infty_0 \mathrm{d}\lambda
\mathrm{Ai}(y+\lambda)\mathrm{Ai}(y'+\lambda)\,.
\end{equation}
The associated extended integral kernel is defined through
\begin{equation}\label{4.5}
K_{{\cal A}_2}(y,\tau;y',\tau')=-(e^{-(\tau'-\tau)H_{\mathrm{Ai}}})(y,y')\mathbf{1}(\tau'>\tau)+(e^{\tau H_{\mathrm{Ai}}} K_{\mathrm{Ai}} e^{-\tau' H_{\mathrm{Ai}}})(y,y').
\end{equation}
Then, the $m$-th marginal for $\mathcal{A}_2(t)$ at times $t_1<t_2<\ldots<t_m$ is expressed as a
determinant on $L^2(\mathbb{R}\times\{t_1,\ldots,t_m\})$ according to~\cite{PS02b}
\begin{equation}\label{4.6}
\mathbb{P}\big(\mathcal{A}_2(t_1)\leq s_1,\ldots,
\mathcal{A}_2(t_m)\leq s_m\big)=\det (\mathbf{1}-\chi_{s}K_{{\cal A}_2}\chi_{s}\big),
\end{equation}
with $\chi_s(x,t_i)=\mathbf{1}(x>s_i)$. $t\mapsto \mathcal{A}_2(t)$ is a stationary process with continuous
sample paths and covariance
$g_2(t)=\mathrm{Cov}\big(\mathcal{A}_2(0),\mathcal{A}_2(t)\big)=\textrm{Var}(\xi_2) - |t|+{\cal O}(t^2)$ for
$t\to 0$ and $g_2(t)=t^{-2}+{\cal O}(t^{-4})$ for $|t|\to\infty$.

The scaling limit for $\mathcal{A}_2(t)$ can be most easily
constructed by two slightly different procedures. The first one
starts from the stationary Ornstein-Uhlenbeck process
$A^{\mathrm{OU}}(t)$ in (\ref{4.2}), which has the potential
$V(x)=x^2/2N$. Its distribution at a single time is GUE, to say
$Z^{-1}_N \exp [-\frac{1}{2N}\mathrm{tr}A^2] \mathrm{d}A$, where the
factor $1/N$ results from the condition that the eigenvalue density in
the bulk is of order~$1$. Let $\lambda^{\mathrm{OU}}_N(t)$ be the
largest eigenvalue of $A^{\mathrm{OU}}(t)$. Then
\begin{equation}\label{4.8}
\lim_{N\to\infty} N^{-1/3} \big(\lambda^{\mathrm{OU}}_N
(2N^{2/3}t)-2N\big)=\mathcal{A}_2(t)
\end{equation}
in the sense of convergence of finite dimensional distributions. The
$N^{2/3}$ scaling means that locally $\lambda^{\mathrm{OU}}_N(t)$
looks like Brownian motion. On the other side the global behavior is confined.

The marginal of the stationary Ornstein-Uhlenbeck process for two
time instants is the familiar 2-matrix model~\cite{EM98,NF98}. Setting
$A_1=A^{\mathrm{OU}}(0)$, $A_2=A^{\mathrm{OU}}(t)$, $t>0$, the joint
distribution is given by
\begin{equation}\label{4.9}
\frac{1}{Z^2_N} \exp \Big(-\frac{1}{2N(1-q^2)}\mathrm{tr}[A^2_1+A^2_2-2qA_1A_2]\Big)
\mathrm{d}A_1 \mathrm{d} A_2,\quad q=\exp(-t/2N).
\end{equation}

A somewhat different construction uses the Brownian bridge defined
by (\ref{4.2}) with $V=0$ and
$A^{\mathrm{BB}}(-T)=A^{\mathrm{BB}}(T)=0$, that is
\begin{equation}\label{4.10}
A^{\mathrm{BB}}(t)= B(T+t)-\frac{T+t}{2T} B(2T)\,,\quad |t|\leq T\,.
\end{equation}
The eigenvalues
$t\mapsto\big(\lambda^{\mathrm{BB}}_1(t),\ldots,\lambda^\mathrm{BB}_N(t)\big)$
is the Brownian bridge line ensemble. Its largest eigenvalue,
$\lambda^{\mathrm{BB}}_N(t)$, has the scaling limit, for $T=2N$,
\begin{equation}\label{4.11}
\lim_{N\to \infty} N^{-1/3} \big(\lambda_N^{\mathrm{BB}} (2N^{2/3}
t)-2N\big)+t^2= \mathcal{A}_2 (t)
\end{equation}
in the sense of finite-dimensional distributions. Note that
$\lambda^{\mathrm{BB}}_N(t)$ is curved on the macroscopic scale
resulting in the displacement by $-t^2$. But with this subtraction the
limit is stationary.

To prove the limits (\ref{4.8}) and (\ref{4.11}) one uses in a central
way that the underlying line ensemble are determinantal. For the
Brownian bridge this can be seen by the following construction. We
consider $N$ independent standard Brownian bridges over $[-T,T]$,
$b^{\mathrm{BB}}_j(t)$, $j=1,\ldots,N$, $b^{\mathrm{BB}}_j(\pm T)=0$
and condition them on non-crossing for $|t|<T$. The resulting line
ensemble has the same statistics as $\lambda^{\mathrm{BB}}_j(t)$,
$|t|\leq T$, $j=1,\ldots,N$, which hence is determinantal.

The TASEP, and its limits, have also an underlying line ensemble,
which qualitatively resemble $\{\lambda^{\mathrm{BB}}_j(t)$,
$j=1,\ldots,N\}$. The construction of the line ensemble is not
difficult, but somewhat hidden. Because of lack of space we explain
only the line ensemble for the PNG droplet. As before $t$ is the
growth time and $x$ is space which takes the role of $t$ from above.
The top line is $\lambda_0(x,t)=h(x,t)$, $h$ the PNG droplet of
Section~\ref{sect1}. Initially we add the extra lines
$\lambda_j(x,0)=j$, $j=-1,-2,\ldots$. The motion of these lines is
completely determined by $h(x,t)$ through the following simple rules:
(1) and (2) from above are in force for all lines $\lambda_j(x,t)$,
$j=0,-1,\ldots$. (3) holds only for \mbox{$\lambda_0(x,t)=h(x,t)$}. (4) The
annihilation of a pair of an adjacent down-step and up-step at line
$j$ is carried out and copied instantaneously as a nucleation event
to line $j-1$.

We let the dynamics run up to time $t$. The line ensemble at time $t$ is
$\{\lambda_j(x,t),|x|\leq t, j=0,-1,\ldots\}$. Note that
$\lambda_j(\pm t,t)=j$. Also, for a given realization of $h(x,t)$,
there is an index $j_0$ such that for $j<j_0$ it holds
$\lambda_j(x,t)=j$ for all $x, |x|\leq t$. The crucial point of the
construction is to have the statistics of the line ensemble
determinantal, a property shared by the Brownian bridge over
$[-t,t]$, $\lambda^{\mathrm{BB}}_j (x)$, $|x|\leq t$. The multi-line
PNG droplet allows for a construction similar to the Brownian bridge
line ensemble. We consider a family of independent, rate 1,
time-continuous, symmetric nearest neighbor random walks on
$\mathbb{Z}$, $r_j(x)$, $j=0,-1,\ldots$. The $j$-th random walk is
conditioned on $r_j(\pm t)=j$ and the whole family is conditioned on
non-crossing. The resulting line ensemble has the same statistics as
the PNG line ensemble $\{\lambda_j(x,t),j=0,-1,\ldots\}$, which hence
is determinantal.

In the scaling limit for $x=\mathcal{O}(t^{2/3})$ the top line
$\lambda_0(x,t)$ is displaced by $2t$ and $t^{1/3}$ away from
$\lambda_{-1}(x,t)$. Similarly $\lambda^{\mathrm{BB}}_N(x)$ for
$x=\mathcal{O}(N^{2/3})$ is displaced by $2N$ and order $N^{1/3}$ apart from
$\lambda^{\mathrm{BB}}_{N-1}(x)$. On this scale the difference
between random walk and Brownian motion disappears but the
non-crossing constraint persists. Thus it is no longer a surprise
that
\begin{equation}\label{4.12}
\lim_{t\to\infty} t^{-1/3} \big(h(t^{2/3} x,t)-2t\big)+
x^2=\mathcal{A}_2 (x)\,,\quad x\in\mathbb{R}\,,
\end{equation}
in the sense of finite dimensional distributions~\cite{PS02b}.

To summarize: for curved initial data the spatial statistics for
large $t$ is identical to the family of largest eigenvalues in a
GUE multi-matrix model.

\section{Flat initial conditions}\label{sect4}

Given the unexpected connection between the PNG model and
GUE multi-matrices, a natural question is whether such a
correspondence holds also for other symmetry classes of random
matrices. The answer is affirmative, but with unexpected twists. Consider the flat PNG model
with \mbox{$h(x,0)=0$} and nucleation events in the whole upper half plane $\{(x,t),
x\in\mathbb{R}, t\geq 0\}$. The removal of the spatial restriction
of nucleation events leads to the problem of the
longest increasing subsequence of a random permutation with involution~\cite{BR01b,PS00b}. The limit shape
will be flat (straight) and in the limit $t\to\infty$ one obtains
\begin{equation}
h(0,t)=2t+\xi_1 t^{1/3}\,,
\end{equation}
where
\begin{equation}
\mathbb{P}(\xi_1\leq s)=F_1(2^{-2/3}s).
\end{equation}
The distribution function $F_1$ is the Tracy-Widom distribution for
the largest GOE eigenvalue.

As before, we can construct the line ensemble and
ask if the link to Brownian motion on GOE matrices persists.
Firstly, let us compare the line ensembles at fixed position for
flat PNG and at fixed time for GOE Brownian motions. In the large time
(resp.\ matrix dimension) limit, the edges of these point processes
converge to the same object: a Pfaffian point process
with $2\times 2$ kernel~\cite{Fer04}. It seems then plausible to conjecture that also
the two line ensembles have the same scaling limit, i.e., the
surface process for flat PNG and for the largest eigenvalue of GOE
Brownian motions should coincide. Since the covariance for the flat PNG has been
computed exactly in the scaling limit, one can compare with simulation results from GOE multi-matrices. The evidence strongly disfavors the conjecture~\cite{BFP08}.

The process describing the largest eigenvalue of GOE multi-matrices
is still unknown, while the limit process of the flat PNG interface
is known~\cite{BFS08} and called the Airy$_1$ process, ${\cal A}_1$,
\begin{equation}
\lim_{t\to\infty} t^{-1/3} \big(h(t^{2/3}x,t)-2t\big)=2^{1/3}{\cal A}_1(2^{-2/3}x).
\end{equation}
Its $m$-point distribution is given in terms of a Fredholm
determinant of the following kernel. Let
$B(y,y')=\textrm{Ai}(y+y')$, $H_1=-\frac{{\rm d}}{{\rm d}y^2}$.
Then,
\begin{equation}
K_{{\cal A}_1}(y,\tau;y',\tau')=-(e^{-(\tau'-\tau)H_1})(y,y')\mathbf{1}(\tau'>\tau)+(e^{\tau H_1} B e^{-\tau' H_1})(y,y')
\end{equation}
and, as for the Airy$_2$ process, the $m$-th marginal for
$\mathcal{A}_1(t)$ at times \mbox{$t_1<t_2<\ldots<t_m$} is expressed
through a determinant on $L^2(\mathbb{R}\times\{t_1,\ldots,t_m\})$
according to
\begin{equation}\label{4.6b}
\mathbb{P}\big(\mathcal{A}_1(t_1)\leq s_1,\ldots,
\mathcal{A}_1(t_m)\leq s_m\big)=\det (\mathbf{1}-\chi_{s}K_{{\cal A}_1}\chi_{s}\big),
\end{equation}
with $\chi_s(x,t_i)=\mathbf{1}(x>s_i)$~\cite{Sas05,BFPS07}. The Airy$_1$ process is a stationary process with covariance
$g_1(t)=\mathrm{Cov}\big(\mathcal{A}_1(0),\mathcal{A}_1(t)\big)=\textrm{Var}(\xi_1) - |t|+{\cal O}(t^2)$ for
$t\to 0$ and $g_1(t)\to 0$ super-exponentially fast as $|t|\to\infty$.

The Airy$_1$ process is obtained using an approach different from
the PNG line ensemble, but still with an algebraic structure
encountered also in random matrices (in the the GUE-minor process~\cite{JN06,OR06}).
We explain the mathematical structure using the continuous
time TASEP as model, since the formulas are the simplest. For a
while we use the standard TASEP representation in terms of
particles. One starts with a formula by Sch\"utz~\cite{Sch97} for the transition
probability of the TASEP particles from generic positions. Consider
the system of $N$ particles with positions
$x_1(t)>x_2(t)>\ldots>x_N(t)$ and let
\begin{equation}
G(x_1,\ldots,x_N;t)=\mathbb{P}(x_1(t)=x_1,\ldots,x_N(t)=x_N |x_1(0)=y_1,\ldots,x_N(0)=y_N).
\end{equation}
Then
\begin{equation}\label{eq5.1b}
G(x_1,\ldots,x_N;t)=\det\left(F_{i-j}(x_{N+1-i}-y_{N+1-j},t)\right)_{1\leq i,j\leq N}
\end{equation}
with
\begin{equation}
F_n(x,t)=\frac{1}{2\pi\textrm{i}}\oint_{|w|>1}{\rm d}w\frac{e^{t(w-1)}}{w^{x-n+1}(w-1)^n}.
\end{equation}
The function $F_n$ satisfies the relation
\begin{equation}
F_n(x,t)=\sum_{y\geq x}F_{n-1}(y,t)\,.
\end{equation}
The key observation is that (\ref{eq5.1b}) can be written as the
marginal of a determinantal line ensemble, i.e. of a measure which
is a product of determinants~\cite{Sas05}. The ``lines'' are denoted by $x^n_i$
with time index $n$, $1\leq n \leq N$, and space index $i$, $1\leq
i\leq n$. We set $x^n_1=x_n$. Then
\begin{equation}\label{eq5.2}
G(x_1,\ldots,x_N;t) = \!\!\! \sum_{x_i^n,2\leq i\leq n\leq N}
\Big(\prod_{n=1}^{N-1} \det(\phi_n(x_i^{n},x_j^{n+1}))_{i,j=1}^n \Big)
\det(\Psi^N_{N-i}(x_j^N))_{i,j=1}^N
\end{equation}
with $\Psi^N_{N-i}(x)=F_{-i+1}(x-y_{N+1-i},t)$, $\phi_n(x,y)=\mathbf{1}(x>y)$
and $\phi_n(x_{n+1}^{n},y)=1$ (here $x_{n+1}^{n}$ plays the role of
a virtual variable). The line ensembles for the PNG droplet and
GUE-valued Brownian motion have the same determinantal structure.
However in (\ref{eq5.2}) the determinants are of increasing sizes which
requires to introduce the virtual variables $x_{n+1}^n$. However, from
the algebraic perspective the two cases can be treated in a similar
way. As a result, the measure in (\ref{eq5.2}) is determinantal (for
instance for any fixed initial conditions, but not only) and has a
defining kernel dependent on $y_1,\ldots,y_N$. The distribution of
the positions of TASEP particles are then given by a Fredholm
determinant of the kernel. To have flat initial conditions, one sets
$y_i=N-2i$, takes first the $N\to\infty$ limit, and then analyzes the
system in the large time limit to get the Airy$_1$ process
defined above.

A couple of remarks:\\
(1) For general initial conditions (for instance for flat initial conditions),
the measure on $\{x_i^n\}$ is not positive, but some projections, like on
the TASEP particles $\{x_1^n\}$, defines a probability measure.\\
(2) The method can be applied also to step initial conditions ($x_k(0)=-k$, $k=1,2,\ldots$)
and one obtains the Airy$_2$ process. In this case, the measure on $\{x_i^n\}$
is a probability measure.\\
(3) In random matrices a measure which is the product of determinants
of increasing size occurs too, for instance in the GUE-minor process~\cite{JN06,OR06}.

\section{Growth models and last passage percolation}\label{sect5} For
the KPZ equation we carry out the Cole-Hopf transformation
\mbox{$Z(x,t)=\exp(-\lambda h(x,t)/\nu_0)$} with the result
\begin{equation}\label{1}
\frac{\partial}{\partial t}Z(x,t)=-\left(-\nu_0
\frac{\partial^2}{\partial x^2}+\frac{\lambda}{\nu_0}\eta (x,t)\right)
Z(x,t)\, ,
\end{equation}
which is a diffusion equation with random potential.
Using the Feynman-Kac formula, (\ref{1}) corresponds to Brownian
motion paths, $x_t$, with weight
\begin{equation}\label{2}
\exp \left(-\frac{\lambda}{\nu_0}\int^t_0 \mathrm{d}s\, \eta(x_s,s)\right).
\end{equation}
In physics this problem is known as a directed polymer in a random
potential, while in probability theory one uses directed first/last passage
percolation. The spirit of the somewhat formal expression (\ref{2}) persists for
discrete growth processes. For example, for the PNG droplet we fix a
realization, $\omega$, of the nucleation events, which then
determines the height $h(0,t;\omega)$ according to the PNG rules. We
now draw a piecewise linear path with local slope $m$, $|m|<1$,
which starts at $(0,0)$, ends at $(0,t)$, and changes direction only
at the points of $\omega$. Let $L(t;\omega)$ be the maximal number
of Poisson points collected when varying over allowed paths. Then
$h(0,t;\omega)=L(t;\omega)$. So to speak, the random potential from (\ref{2}) is
spiked at the Poisson points.

In this section we explain the connection between growth models and
(directed) last passage percolation. For simplicity, we first
consider the case leading to discrete time PNG droplet, although
directed percolation can be defined for general passage time distributions. Other
geometries like flat growth are discussed later.

Let $\omega(i,j)$, $i,j\geq 1$, be independent random variables with
geometric distribution of parameter $q$, i.e.,
$\mathbb{P}(\omega(i,j)=k)=(1-q)q^k$, $k\geq 0$. An up-right path
$\pi$ from $(1,1)$ to $(n,m)$ is a sequence of points
$(i_\ell,j_\ell)_{\ell=1}^{n+m-1}$ with
$(i_{\ell+1},j_{\ell+1})-(i_\ell,j_\ell)\in \{(1,0),(0,1)\}$. The
last passage time from $(1,1)$ to $(n,m)$ is defined by
\begin{equation}\label{eq5.1}
G(n,m)=\max_{\pi:(1,1)\to (n,m)}\sum_{(i,j)\in \pi} \omega(i,j).
\end{equation}
The connection between directed percolation and different growth
models is obtained by appropriate projections of the
three-dimensional graph of $G$. Let us see how this works.

\begin{figure}
\begin{center}
\epsfig{file=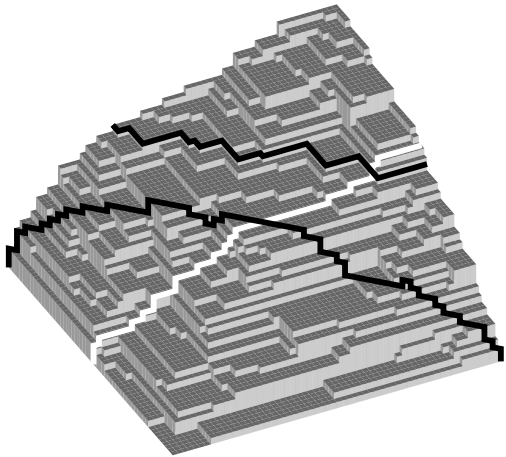,height=6cm}
\caption{Directed percolation and growth models}
\label{FigDirPerc}
\end{center}
\end{figure}

Let time $t$ be defined by $t=i+j-1$ and position by $x=i-j$. Then,
the connection between the height function of the discrete time PNG
and the last passage time $G$ is simply~\cite{Joh03}
\begin{equation}
h(x,t)=G(i,j).
\end{equation}
Thus, the discrete PNG droplet is nothing else than the time-slicing along the $i+j=t$ directions, see Figure~\ref{FigDirPerc}.

We can however use also a different slicing, at constant $\tau=G$, to
obtain the TASEP at time $\tau$ with step initial conditions. For
simplicity, consider $\omega(i,j)$ to be exponentially distributed
with mean one. Then $\omega(n,m)$ is the waiting time for particle
$m$ to do his $n$th jump. Hence, $G(n,m)$ is the time when particle $m$
arrives at $-m+n$, i.e.,
\begin{equation}
\mathbb{P}(G(n,m)\leq \tau)=\mathbb{P}(x_m(\tau)\geq -m+n).
\end{equation}

From the point of view of the TASEP, there is another interesting
cut, namely at fixed $j=n$. This corresponds to look at the
evolution of the position of a given (tagged) particle, $x_n(\tau)$.

\newpage
A few observations:\\
(1) Geometrically distributed random passage times correspond to discrete time TASEP
with sequential update.\\
(2) The discrete time TASEP with parallel update is obtained by
replacing $\omega(i,j)$ by $\omega(i,j)+1$.

The link between last passage percolation and growth holds also for
general initial conditions. In (\ref{eq5.1}) the optimization
problem is called point-to-point, since both $(1,1)$ and $(n,m)$ are
fixed. We can generalize the model by allowing $\omega(i,j)$ to be
defined on $(i,j)\in\mathbb{Z}^2$ and not only for $i,j\geq 1$.
Consider the line $L=\{i+j=2\}$ and the following point-to-line
maximization problem:
\begin{equation}
G_L(n,m)=\max_{\pi:L\to (n,m)}\sum_{(i,j)\in \pi} \omega(i,j).
\end{equation}
Then, the relation to the discrete time PNG droplet, namely
\mbox{$h(x,t)=G_L(i,j)$} still holds but this time $h$ is the height
obtained from flat initial conditions. For the TASEP, it means to
have at time $\tau=0$ the particles occupying $2\mathbb{Z}$. Also
random initial conditions fit in the picture, this time one has to
optimize over end-points which are located on a random line.

In the appropriate scaling limit, for large time/particle number one
obtains the Airy$_2$ (resp.\ the Airy$_1$) process for all these
cases. One might wonder why the process seems not to depend on the
chosen cut. In fact, this is not completely true. Indeed, consider
for instance the PNG droplet and ask the question of joint
correlations of $h(x,t)$ in space-time. We have seen that for large
$t$ the correlation length scales as $t^{2/3}$. However, along the
rays leaving from $(x,t)=(0,0)$ the height function decorrelates on
a much longer time scale, of order $t$. These slow decorrelation
directions depend on the initial conditions. For instance, for flat
PNG they are the parallel to the time axis. More generally, they coincide
with the characteristics of the macroscopic surface evolution.
Consequently, except along the slow directions, on the $t^{2/3}$
scale one will always see the Airy processes.

\section{Growth models and random tiling}\label{sect6} In the previous
section we explained how different growth models (TASEP and PNG) and
observables (TASEP at fixed time or tagged particle motion) can be
viewed as different projections of a single three-dimensional
object. A similar unifying approach exists also for some growth
models and random tiling problems: there is a $2+1$ dimensional
surface whose projection to one less dimension in space (resp.\
time) leads to growth in 1+1 dimensions (resp.\ random tiling in 2
dimensions)~\cite{BF08a}. To explain the idea, we consider the dynamical model
connected to the continuous time TASEP with step initial conditions,
being one of the simplest to define.

\begin{figure}
\begin{center}
\psfrag{x}[b]{$x$}
\psfrag{n}[b]{$n$}
\psfrag{h}[b]{$h$}
\psfrag{x11}[][][0.9]{$x_1^1$}
\psfrag{x12}[][][0.9]{$x_1^2$}
\psfrag{x13}[][][0.9]{$x_1^3$}
\psfrag{x14}[][][0.9]{$x_1^4$}
\psfrag{x22}[][][0.9]{$x_2^2$}
\psfrag{x23}[][][0.9]{$x_2^3$}
\psfrag{x24}[][][0.9]{$x_2^4$}
\psfrag{x33}[][][0.9]{$x_3^3$}
\psfrag{x34}[][][0.9]{$x_3^4$}
\psfrag{x44}[][][0.9]{$x_4^4$}
\psfrag{tasep}[c]{\textbf{TASEP}}
\epsfig{file=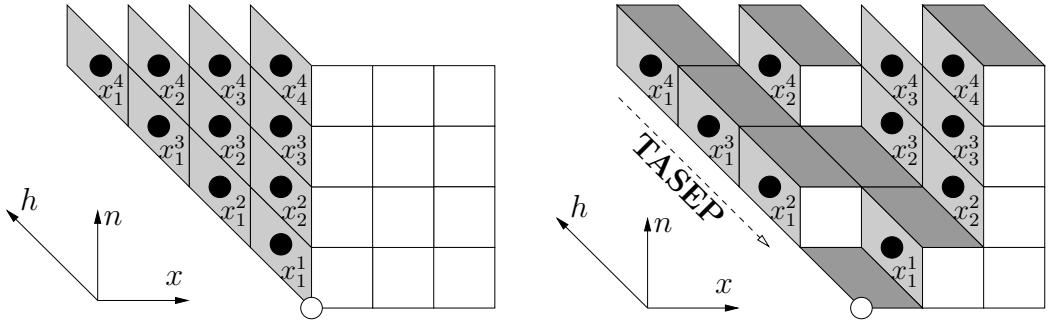,width=\textwidth} \caption{Illustration
of the particle system and its corresponding lozenge tiling. The initial condition is on the left.}
\label{FigDynamics2d}
\end{center}
\end{figure}

In Section~\ref{sect4} we encountered a measure on set of variables
\mbox{$S_N=\{x_i^n,1\leq i\leq n\leq N\}$}, see (\ref{eq5.2}). The product
of determinants of the $\phi_n$'s implies that the measure is
non-zero, if the variables $x_i^n$ belong to the set $S_N^{\rm int}$
defined through an interlacing condition,
\begin{equation}
S_N^{\rm int}=\{x_i^n \in S_N \,|\, x_i^{n+1}<x_i^n\leq x_{i+1}^{n+1} \}.
\end{equation}
Moreover, for TASEP with step initial conditions, the measure on
$S_N^{\rm int}$ is a probability measure, so that we can interpret
the variable $x_i^n$ as the position of the particle indexed by
$(i,n)$. Also, the step initial condition, \mbox{$x_n(t=0)\equiv x_1^n(0)=-n$}, implies that $x_i^n(0)=i-n-1$, see
Figure~\ref{FigDynamics2d} (left).

Then, the dynamics of the TASEP induces a dynamics on the particles
in $S_N^{\rm int}$ as follows. Particles independently try to jump
on their right with rate one, but they might be blocked or pushed by
others.
When particle $x_k^n$ attempts to jump:\\[0.5em]
(1) it jumps if $x_i^n<x_i^{n-1}$, otherwise it is blocked
(the usual TASEP dynamics between particles with same lower index),\\[0.5em]
(2) if it jumps, it pushes by one all other particles with index
$(i+\ell,n+\ell)$, $\ell\geq 1$, which are at the same position
(so to remain in the set $S_N^{\rm int}$).\\[0.5em]
For example, consider the particles of Figure~\ref{FigDynamics2d}
(right). Then, if in this state of the system particle $(1,3)$ tries
to jump, it is blocked by the particle $(1,2)$, while if particle
$(2,2)$ jumps, then also $(3,3)$ and $(4,4)$ will move by one unit
at the same time.

It is clear that the projection of the particle system onto
\mbox{$\{x_1^n,1\leq n\leq N\}$} reduces to the TASEP dynamics in
continuous time and step initial conditions, this projection being
the sum in (\ref{eq5.2}).

\begin{figure}
\begin{center}
\psfrag{z11}{$z_1^1$}
\psfrag{z12}{$z_1^2$}
\psfrag{z22}{$z_2^2$}
\psfrag{t0}[c]{$t=0$}
\psfrag{t1}[c]{$t=1$}
\psfrag{t2}[c]{$t=2$}
\epsfig{file=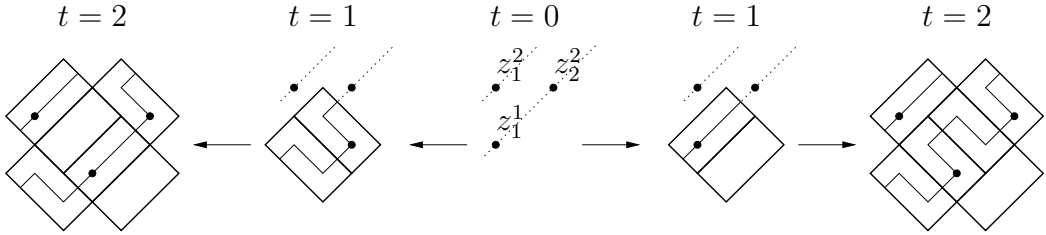,width=\textwidth}
\caption{Two examples of configurations at time $t=2$
obtained by the particle dynamics and its associated Aztec
diamonds. The line ensembles are also drawn.}
\label{FigAztec}
\end{center}
\end{figure}

The system of particles can also be represented as a tiling using
three different lozenges as indicated in Figure~\ref{FigDynamics2d}.
The initial condition corresponds to a perfectly ordered tiling and
the dynamics on particles reflects a corresponding dynamics of the
random tiling. Thus, the projection of the model to fixed time reduces to a random tiling problem.
In the same spirit, the shuffling algorithm of the Aztec diamond
falls into place. This time the interlacing is $S_{\rm Aztec}=\{z_i^n \,|\,
z_i^{n+1}\leq z_i^n\leq z_{i+1}^{n+1}\}$ and the dynamics is on
discrete time, with particles with index $n$ not allowed to move
before time $t=n$. As before, particle $(i,n)$ can be blocked by
$(i,n-1)$. The pushing occurs in the following way: particle $(i,n)$ is forced to move
if it stays at the same position of $(i-1,n-1)$, i.e., if it would generate a violation
due to the possible jump of particle $(i-1,n-1)$. Then at
time $t$ all particles which are not blocked or forced to move, jump
independently with probability $q$. As explained in
detail in~\cite{Nor08} this particle dynamics is equivalent to the
shuffling algorithm of the Aztec diamond (take $q=1/2$ for uniform
weight). In Figure~\ref{FigAztec} we illustrate the rules with a small size
example of two steps. There we also draw a set of lines,
which come from a non-intersecting line ensemble similar to the ones of the
PNG model and the matrix-valued Brownian motions described in
Section~\ref{sect3}.

On the other hand, let $x_i^n:=z_i^n-n$. Then, the dynamics of the
shuffling algorithm projected onto particles $x_1^n$, $n\geq 1$, is
nothing else than the discrete time TASEP with parallel update and
step initial condition! Once again, we have a $1+1$ dimensional growth and a
$2$ dimensional tiling model which are different projections of the same
$2+1$-dimensional dynamical object.

\section{A guide to the literature}\label{sect7}

There is a substantial body of literature and only a rather partial list is given here. The guideline is ordered according to physical model under study.\smallskip\\

\textit{PNG model.} A wide variety of growth processes, including
PNG, are explained in~\cite{Mea98}. The direct link to the largest
increasing subsequences of a random permutation and to the unitary
matrix integral of~\cite{PS90} is noted in~\cite{PS00b,PS00a}. The
convergence to the Airy$_2$ process is worked out in~\cite{PS02b} and
the stationary case in~\cite{BR00,PS04}. For flat initial
conditions, the height at a single space-point is studied
in~\cite{BR01a,BR01b} and for the ensemble of top lines
in~\cite{Fer04}. The extension to many space-points is accomplished
by~\cite{BFS08}. Determinantal space-time processes for the discrete time PNG model are discussed
by~\cite{Joh03}. External source at the origin is studied in~\cite{IS04b} for the full line and in~\cite{IS04a} for the half-line.
\\

\textit{Asymmetric simple exclusion }(ASEP). As a physical model reversible lattice gases, in particular the simple exclusion process, were introduced by Kawasaki~\cite{Kaw72} and its asymmetric version
by Spitzer~\cite{Spi70}. We refer to Liggett~\cite{Lig99} for a
comprehensive coverage from the probabilistic side. The hydrodynamic
limit is treated in~\cite{Spo91} and~\cite{KL99}, for example. There
is a very substantial body on large deviations with Derrida as one of the
driving forces, see~\cite{Der07} for a recent review. For the TASEP Sch\"utz~\cite{Sch97} discovered a determinantal-like formula for the transition probability for a finite number of particles on $\mathbb{Z}$. TASEP step initial conditions are studied in the seminal paper by
Johansson~\cite{Joh00}. The random matrix representation of~\cite{Joh00} may also be obtained by the Sch\"utz formula~\cite{NS04,RS05}. The convergence to the
Airy$_2$ process in a discrete time setting~\cite{Joh05}. General step
initial conditions are covered by~\cite{PS02a} and the extended process in~\cite{BFP09}. In~\cite{FS06b} the
scaling limit of the stationary two-point function is proved.
Periodic intial conditions were first studied by
Sasamoto~\cite{Sas05} and widely extended in~\cite{BFPS07}. A
further approach comes from the Bethe ansatz which is used
in~\cite{GS92} to determine the spectral gap of the generator.
In a spectacular analytic tour de force Tracy and Widom develop the Bethe ansatz
for the transition probability and thereby extend the Johansson result to the
PASEP~\cite{TW09,TW08}.\\

\textit{2D tiling (statics).} The connection between growth and
tiling was first understood in the context of the Aztec
diamond~\cite{JPS98}, who prove the arctic circle theorem. Because of the
specific boundary conditions imposed, for typical
tilings there is an ordered zone spatially coexisting with a
disordered zone. In the disordered zone the large scale statistics
is expected to be governed by a massless Gaussian field, while the
line separating the two zones has the Airy$_2$ process
statistics.\medskip\\
a) Aztec diamond. The zone boundary is studied by~\cite{JPS98} and by~\cite{JN06}.
Local dimer statistics are investigated in~\cite{CEP96}. Refined details are the large scale
Gaussian statistics in the disordered zone~\cite{Ken00}, the edge statistics~\cite{Joh05},
and the statistics close to a point touching the boundary~\cite{JN06}.\medskip\\
b) Ising corner. The Ising corner corresponds to a lozenge tiling under the constraint of a fixed volume below the thereby defined surface. The largest scale information is the macroscopic shape and large deviations~\cite{CK01}. The determinantal structure is noted in~\cite{OR03}. This can be used to study the edge statistics~\cite{FS03,FPS04}. More general boundary conditions (skew plane partitions) leads to a wide variety of macroscopic shapes~\cite{OR07}.\medskip\\
c) Six-vertex model with domain wall boundary conditions, as
introduced in~\cite{KZJ00}. The free energy including prefactors is available~\cite{BF06}.
A numerical study can be found in~\cite{AR05}. The mapping to the
Aztec diamond, on the free Fermion line, is explained
in~\cite{FS06a}.\medskip\\
d) Kasteleyn domino tilings. Kasteleyn~\cite{Kas63} noted that
Pfaffian methods work for a general class of lattices. Macroscopic
shapes are obtained by~\cite{KOS06} with surprising connections to
algebraic geometry. Gaussian fluctuations are proved
in~\cite{Ken08}.\\

\textit{2D tiling (dynamics)}, see Section~\ref{sect6}. The shuffling
algorithm for the Aztec diamond is studied
in~\cite{EKLP92,Pro03,Nor08}. The pushASEP is introduced
by~\cite{BF08b} and anisotropic growth models are investigated
in~\cite{BF08a}, see also~\cite{PS97} for the Gates-Westcott model. A
similar intertwining structure appears for Dyson's Brownian
motions~\cite{War07}.\\

\textit{Directed last passage percolation}. ``Directed'' refers to
the constraint of not being allowed to turn back in time. In the physics
literature one speaks of a directed polymer in a random medium.
Shape theorems are proved, e.g., in~\cite{Kes86}. While the issue of
fluctuations had been repeatedly raised, sharp results had to wait
for~\cite{Joh00} and~\cite{PS02b}. Growths models naturally lead to
either point-to-point, point-to-line, and point-to-random-line
last passage percolation. For certain models one has to further impose
boundary conditions and/or symmetry
conditions for the allowed domain. In (\ref{eq5.1}) one takes
the max, thus zero temperature. There are interesting results for
the finite temperature version, where the energy appear in the exponential, as in (\ref{2})~\cite{CY06}.\\

\textit{KPZ equation.} The seminal paper is~\cite{KPZ86}, which
generated a large body of theoretical work. An introductory
exposition is~\cite{BS95}. The KPZ equation is a stochastic field
theory with broken time reversal invariance, hence a great
theoretical challenge, see, e.g.,~\cite{Las98}.\\

\textit{Review articles.} A beautiful review is~\cite{Joh06}. Growth
models, of the type discussed here, are explained in~\cite{FP06,Fer08}.
A fine introduction to random matrix techniques
is~\cite{Sas07}.~\cite{KK08} provide an introductory exposition to
the shape fluctuation proof of Johansson. The method of line
ensembles is reviewed in~\cite{Spo06}.

\end{document}